\documentclass{amsart}
\usepackage{graphicx}
\theoremstyle{definition}

\theoremstyle{remark}

\numberwithin{equation}{section}

\begin{document}

\title{A remark on odd dimensional normalized Ricci flow}%
\author{Hong Huang}%
\address{Department of Mathematics,Beijing Normal University,Beijing 100875, P. R. China}%
\email{hhuang@bnu.edu.cn}%

\thanks{Partially supported by NSFC no.10671018.}%

\begin{abstract}
Let $(M^n,g_0)$ ($n$ odd) be a compact Riemannian manifold with
$\lambda(g_0)>0$, where $\lambda(g_0)$ is the first eigenvalue of
the operator $-4\Delta_{g_0}+R(g_0)$, and $R(g_0)$ is the scalar
curvature of $(M^n,g_0)$. Assume the maximal solution $g(t)$ to the
normalized Ricci flow with initial data $(M^n,g_0)$ satisfies
$|R(g(t))| \leq C$ and $\int_M |Rm(g(t))|^{n/2}d\mu_t \leq C$
uniformly for a constant $C$. Then we show that the solution
sub-converges to a shrinking Ricci soliton. Moreover,when $n=3$, the
condition $\int_M |Rm(g(t))|^{n/2}d\mu_t \leq C$ can be removed.
\end{abstract} \maketitle

Since Hamilton's seminal work [H1] Ricci flow has been an important
tool used extensively in  geometry and topology.  In particular,
there is the recent breakthrough of Perelman [P1],[P2].

In this short note we  prove a convergence result for odd
dimensional volume-normalized Ricci flow,

$\frac{\partial g_{ij}}{\partial t}=-2R_{ij}+\frac{2}{n}rg_{ij}$,

where $r=\frac{\int_M Rd\mu_t}{\int_M d\mu_t}$ is the average
scalar curvature of $(M^n,g(t))$. More precisely, we have the
following

\hspace *{0.4cm}

{\bf Theorem } Let $(M^n,g_0)$ ($n$ odd) be a compact Riemannian
manifold with $\lambda(g_0)>0$, where $\lambda(g_0)$ is the first
eigenvalue of the operator $-4\Delta_{g_0}+R(g_0)$. Assume the
maximal solution $g(t)$ to the volume-normalized Ricci flow with
initial data $(M^n,g_0)$ satisfies $|R(g(t))| \leq C$ and $\int_M
|Rm(g(t))|^{n/2}d\mu_t \leq C$ uniformly  for a constant $C$. Then
the solution sub-converges to a shrinking Ricci soliton.
Moreover,when $n=3$, the condition $\int_M |Rm(g(t))|^{n/2}d\mu_t
\leq C$ can be removed.

 (Here, $R(g_0)$ is the scalar curvature of $(M^n,g_0)$,
and $Rm(g(t))$ is the curvature tensor of $g(t)$.)

\hspace *{0.4cm}

{\bf Proof} Let $[0,T)$ be the maximal time interval of existence
of $g(t)$. First we show that $|Rm(g(t))|\leq C'$ uniformly on
$[0,T)$ for a constant $C'$. Suppose this is not the case, then
there exist a sequence of times $t_i\rightarrow T$ and points
$x_i\in M$ such that
      $Q_i=|Rm(g(t_i))|(x_i)=max_{x\in M} |Rm(g(t_i))|(x)\rightarrow \infty$.
  Note that the condition $\lambda(g_0)>0$ implies the corresponding un-normalized solution
  blows up in finite time, since we have $\frac{d}{d\tilde{t}}\lambda(\tilde{t}) \ge \frac{2}{n}\lambda^{2}(\tilde{t})$
  by Perelman [P1], where $\lambda(\tilde{t})$ is the first eigenvalue of the operator $-4\Delta+R$ for the un-normalized
  Ricci flow. By  Perelman's no local collapsing theorem [P1] and
  Hamilton's compactness theorem [H2], a subsequence
  of the rescaled solutions $(M,Q_ig({Q_i}^{-1}t+t_i),x_i)$ converges smoothly to a
  pointed complete (``normalized"-)Ricci flow $(M_\infty,g_\infty(t),x_\infty)$, such that
   $g_\infty(t)$ is $\kappa$-noncollapsed relative to upper bounds of the scalar curvature
   on all scales, where $\kappa$ is certain positive constant depending only on $n$ and
    the initial data $g_0$.   Clearly $M_\infty$ is non-compact.
    We have
    $R(g_\infty(t))=0$, hence $g_\infty(t)$ is
    Ricci flat since it is a solution to the Ricci
    flow.
    Moreover, the conditions (ii) and (iii) in (3.14) of [A] are also
    fulfilled for $(M_\infty, g_\infty(0))$. Combined with the odd dimensional assumption, it
    follows from [A, Theorem 3.5](see also [BKN,Theorem(1.5)]) that (the double cover of )
    $(M_\infty, g_\infty(0))$ is the $n$-dimensional Euclidean space, which
    contradicts  the fact $|Rm(g_\infty(0))|(x_\infty)=1$.

When $n=3$ and the condition $\int_M |Rm(g(t))|^{n/2}d\mu_t \leq C$
is removed, then if the curvature tensors of $g(t)$ are not bounded
uniformly on $[0,T)$, the $(M_\infty, g_\infty(0))$ constructed
above is Ricci flat and hence flat, and again we get a
contradiction.

Then it follows that $T=\infty$ and  $g(t)$ is nonsingular. By [FZZ,
Proposition 2.2] for any sequence of times $t_k \rightarrow \infty$,
there is a subsequence $t_{k_i}$ such that $g(t+t_{k_i})$ converges
to a shrinking Ricci soliton. Moreover,when $n=3$, we must have that
$M$ is diffeomorphic to a spherical space form (cf. [H3]).

\hspace *{0.4cm}

{\bf Remark 1 } A similar argument was used by Ruan, Zhang and Zhang
in [RZZ], where they proved a related result (see Proposition 1.3 in
[RZZ]).

{\bf Remark 2 } Of course, compact three manifolds $(M^3,g_0)$ with
$\lambda(g_0)>0$ have been classified by Perelman [P2].

\hspace *{0.4cm}

{

\bibliographystyle{amsplain}

\hspace *{0.4cm}

{\bf Reference}

\bibliography{1}[A] M.T.Anderson,Ricci curvature bounds and Einstein metrics
 on compact manifolds, J. Amer. Math. Soc. 2(1989),455-490.

\bibliography{2}[BKN] S. Bando, A. Kasue, and H. Nakajima, On a construction of coordinates
at infinity on manifolds with fast curvature decay and maximal
volume growth, Invent. Math. 97(1989), 313-349.

\bibliography{3}[FZZ] F. Fang, Y. Zhang, and Z. Zhang, Non-singular solutions
to the normalized Ricci flow equation, arXiv:math.GM/0609254.

\bibliography{4}[H1] R. S. Hamilton,Three-manifolds with positive Ricci curvature,
J. Diff. Geom. 17 (1982),255-306.

\bibliography{5}[H2] R. S. Hamilton, A compactness property for solutions of
 the Ricci flow, Amer. J. Math. 117(1995), 545-572.

\bibliography{6}[H3] R. S. Hamilton, Non-singular solutions of the Ricci flow on three-manifolds,
Commun. Anal. Geom. 7(1999), 695-729.

\bibliography{7}[P1] G. Perelman, The entropy formala for the Ricci
flow and its geometric applications, arXiv:math.DG/0211159.

\bibliography{8}[P2] G. Perelman, Ricci flow with surgery on three
manifolds, arXiv:math.DG/0303109.

\bibliography{9}[RZZ]W.-D. Ruan, Y. Zhang, and Z. Zhang, Bounding sectional curvature along
a K$\ddot{a}$hler-Ricci flow, arXiv:0710.3919.

\end{document}